\documentclass[twocolumn]{elsart}
\usepackage{amssymb,amsmath,graphicx}

\makeatletter
\renewcommand{\@makefnmark}{\hbox{\mathsurround=0pt}$^{\mbox{\dag}}$}
\makeatother

\begin{document}

\textbf{\LARGE One-dimensional chaos in a system with dry \\  friction: analytical approach.}

\bigskip

{\noindent\large Nikita Begun${}^a$, Sergey Kryzhevich${}^{a,b,c}$\footnote{Corresponding author.\\\hspace*{5.5mm}Email address:
kryzhevicz@gmail.com, phone: +79219181076, fax: +78124286944}~}

\medskip

{\noindent\small
${}^a$ Faculty of Mathematics and Mechanics, Saint-Petersburg State University,\\
28, Universitetskiy pr., Peterhof, Saint-Petersburg, 198503, Russia;\\
${}^b$ University of Aveiro, Department of Mathematics, 3810$-$193, Aveiro, Portugal;\\
${}^c$ BCAM -- Basque Center for Applied Mathematics, Mazarredo, 14, E48009 Bilbao, Basque-Country, Spain.}

\bigskip

\noindent\textbf{Abstract.} We introduce a new analytical method, which allows to find out chaotic dynamics in non-smooth dynamical systems. A simple mechanical system consisting of a mass and a dry friction element is considered as an example. The corresponding mathematical model is represented. We show that the considered dynamical system is a skew product over a piecewise smooth mapping of a segment (a base map). For this base map we demonstrate there is a domain of parameters where a robust chaotic dynamics can be observed. Namely, we prove existence of an infinite set of periodic points of arbitrarily big period. Moreover, a reduction of the considered map is semi-conjugated to a shift on the set of one-sided infinite boolean sequences. Also, we find conditions, sufficient for existence of a superstable periodic point of this map. The obtained result partially solves a general problem: theoretical confirmation of chaotic and periodic regimes numerically and experimentally observed for models of percussion drilling.

\bigskip


\bigskip


\bigskip

\noindent\textbf{1. Introduction.}

\medskip

Systems with dry friction form a wide and important class inside discontinuous dynamical systems. They appear in many applications, especially in manufacturing systems: vibrating conveyors, percussion drilling, metal cutting, etc (see \cite{kb1,kb2,kb3,kb4,kb5,kb6,kb7,kb8,kb9,kb10,kb11,kb12} and references therein for review). Their properties manifest many principle differences with ones of smooth dynamics. For instance, the uniqueness theorem is not valid any more.  An approach to study such systems has been developed by A.\,F.\, Filippov \cite{kb13}. He offered to consider piecewise continuous systems of differential equations as families of vector fields, defined on disjoint domains of the phase space and define auxiliary tangent flows on boundaries, respecting limit directions of vector fields. This approach reduces a discontinuous system to differential inclusions. Moreover, the phase space may become multidimensional e.g. a set of initial data of the full dimension may be transferred to a set of a lower dimension. The theory of discontinuous systems and specific bifurcation is well-developed \cite{kb1,kb2,kb4,kb8,kb9,kb11,kb13,kb14,kb15}. It is also well-known that chaotic dynamics frequently occurs in such systems, particularly in ones with dry friction \cite{kb1,kb2,kb3,kb5,kb9,kb10,kb11,kb12,kb16,kb17,kb18,kb19}.

Apart from numerical and experimental simulations, the most common analytic approach involves a reduction of dimension. For some systems with dry friction it is possible to demonstrate that there exists an invariant set of dimension 1 where the attractor resides. A method to find this attractor has been developed by M.\,Wiercigroch, E.\,Pavlovskaya and A.Krivtsov in papers \cite{kb6,kb7} and, in its general form, in the paper \cite{kb18}. In our paper we use some ideas of this approach. Another powerful method has been proposed by R.\,Szalai and H.\,M.\,Osinga \cite{kb19}. They have proved that for a general class of systems with dry friction the attractor resides in a polygon type set and demonstrated a possibility of a chaotic dynamics there. Later \cite{kb20} they have shown, using a modification of their method that some complex structures like Arnold tongues can be observed in a neighborhood of the so-called grazing-sliding bifurcation \cite{kb1,kb2,kb15}.

The main aim of this paper is to provide a new method which allows to find chaotic invariant sets in systems with a dry friction. To demonstrate this method, we use a very simple example of a system with dry friction, first considered be A.\, Krivtsov and M.\, Wiercigroch \cite{kb6}. First of all, we show that the considered dynamical system engenders a discontinuous mapping of a segment. Here we use ideas from \cite{kb18}. Then we study properties of this mapping which allow us to find two disjoint segments such that the image of every one of them covers their union. Moreover, we prove that the considered mapping is continuous on the union of these segments. This allows us to apply well-known techniques of one dimensional dynamics \cite[Part 3, Section 15]{kb21} and to demonstrate that a kind of chaotic dynamics, similar to one described by T.\,Li and J.\, Yorke \cite{kb22} is there.

The main advantages of the offered method are the following.
\begin{enumerate}
\item We can obtain chaotic sets which, in general, are not attractors.
\item For simple systems with dry fiction, the offered method gives coefficient type criteria of chaos.
\item Though we need a presence of a small parameter in our proofs, it is possible to estimate numerically how small this parameter must be. In general, presence of chaotic invariant sets does not correspond to a neighborhood of any bifurcation.
\item The fact of presence of chaotic behavior (but not corresponding invariant sets) is robust.
\item A corresponding invariant measure can be described using techniques of \cite[Part 3, Section 15]{kb21}.
\end{enumerate}

The paper is organized as follows. In Section 2 we introduce the mathematical model of the considered system and describe possible regimes of motion. In Section 3 we define the main object of our investigation: the one-dimensional mapping, corresponding to phases of switching for solutions. For this, we describe all possible scenarios of behavior of solutions. In the next section we describe some properties of the introduced mapping. In Section 5 we study how the segments of continuity of the constructed mapping look like and find out two segments of continuity whose images cover their union. Sections 6 and 7 are technical. We prove existence of periodic points of all possible periods but we cannot apply the theory of mappings of a segment directly since we deal with a discontinuous mapping. However, still we can use standard methods of this theory to finish our proof. Obtained chaotic invariant sets may coexist with superstable fixed points of the considered map. This is discussed in Section 8. Main results of the paper are formulated in Section 9. Discussion, including notes on the robustness of the obtained set and some plans on the future research is given at Section 10.

\medskip

\noindent\textbf{2. Description of the mathematical model.}

\medskip

Consider a single degree-of-freedom mechanical system, consisting of a point mass and a delimiter with dry friction (Fig.\, 1) which gives a simple model of percussion drilling. This system consists of a unit mass, whose motion is controlled by a harmonic external force $F(t)$ which is a sum of a positive constant component equal to $2b$ and a harmonic component of a positive amplitude $a$ and a period equal to $2\pi$. Also, the considered system includes a delimiter which provides an additional dry friction as soon as the mass reaches it. The maximal value of this dry friction force is $q$. Here and later we always suppose that all considered parameters are non-dimensional.

Our main aim is to prove that provided some additional conditions are satisfied, the dynamics of the considered system is robustly chaotic in topological sense.
Let $x$ be the current position of the mass and $y$ be one of the delimiter. We assume that the inequality $x\le y$ is always satisfied i.e. the mass cannot penetrate through the delimiter and that the delimiter cannot move to the left so the value $y$ is always non-decreasing.

\begin{figure}[!ht]\begin{center}
\includegraphics*[width=3in]{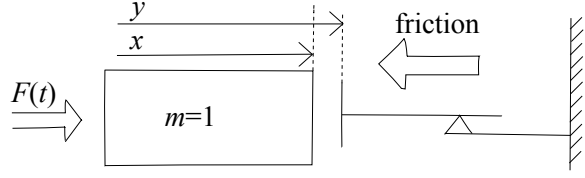}
\end{center}
\caption{The considered mechanical system.}
\end{figure}

Consider the value $\vartheta_0\in [0,\pi/2]$ such that
$$\pi-\vartheta_0=\cot (\vartheta_0/2). \eqno (1)$$
This value is unique and $\vartheta_0\approx 0.81047$, $\sin \vartheta_0 \approx 0.724611$.

We make the following assumptions on parameters of the system:
$$a>0, \qquad b\in (0,a/2), \qquad q\in (a \sin\vartheta_0, a).\eqno (2)$$
We shall always suppose that $b\ll a$ (which implies $b\ll q$). This means that we always suppose that the ratio $b/a$ is as small as necessary.

There are five types of motions of the considered system.

\begin{enumerate}
\item\textbf{No contact (free) motion (f).} This motion takes place if $x<y$ i.e. the mass and the delimiter do not interact. Then
$$\ddot x=F(t)=a\sin t+2b; \quad \dot y=0.\eqno (3)$$
\item\textbf{Contact with progression (p).} In this regime we have an additional friction. The motion is defined by equations
$$\ddot x=F(t)-q=a\sin t+2c=a\sin t+2b-q; \quad y=x.\eqno (4)$$
\item\textbf{Stop (s).} Here the mass and the delimiter are both immobile, i.e.
$$x=y, \qquad \dot x = \dot y=0.\eqno (5)$$
\item\textbf{Instantaneous stop (is).} This happens if Condition (5) is satisfied for a fixed instant of time but is not true in its small neighborhood. So, this happens if the system switches from or to free motion or from/to motion with progression.
\item\textbf{Instantaneous transition from the no contact regime to motion with progression (fp).}
\end{enumerate}

We always suppose that free motion, motion with progression and stop regime are observed at open intervals of time. This allows us to classify all instants of transition.

Let $t_0<t_1$ be zeros of the function $F(t)-q$ and $t_2<t_3$ be ones of the function $F(t)$ on $[0,2\pi]$. Later on we consider the phase
$\varphi=t \mod 2\pi$. Here $\varphi\in S^1={\mathbb R}/2\pi{\mathbb Z}$.

Solutions of Eqs. (3), (4) and (5) can easily be written down. If $x(\theta_0)=x_0$, $\dot x(\theta_0)=x_1$ we have
$$x(t)=-a\sin t+b(t-\theta_0)^2+(x_1+a\cos \theta_0)(t-\theta_0)+x_0+a\sin \theta_0$$
for Eq.\,(3) (free motion) and
$$x(t)=-a\sin t+c(t-\theta_0)^2+(x_1+a\cos \theta_0)(t-\theta_0)+x_0+a\sin \theta_0\eqno (6)$$
for motion with progression. For stop regime we always have $x_1=0$ and $x(t)\equiv x_0$.

\medskip

\noindent\textbf{3. Reduction to dimension 1.}

\medskip
Now we describe how we may proceed from one regime to another.

Starting from the free motion, the mass will always return to the delimiter since $b>0$. If the velocity of the collision is non-zero, the system proceeds to progression regime. The transition $\mathbf{(f)}\to \mathbf{(is)} \to \mathbf{(s)}$ is possible at $t=\theta_0$ if and only if $x(\theta_0)=y(\theta_0)$, $\dot x(\theta_0)=0$. If $F(\theta_0)>0$ then $\dot x(t)$ is negative before $t=\theta_0$ which is impossible. If $F(\theta_0)<0$, the mass returns back, "ignoring"\ presence of the delimiter. We have an instantaneous stop there. Otherwise, $F(\theta_0)=0$. If $\dot F(\theta_0)<0$ then again we have $\dot x(t)<0$ in a left neighborhood of $\theta_0$. So, $\dot F(\theta_0)\ge 0$ and, consequently, $\theta_0=t_3 \mod 2\pi$. Then we have an instantaneous stop and the stop regime later on.

In the motion with progression the derivative $\dot x$ vanishes soon or later since $c<0$. If this happens when $\varphi\in [t_2,t_3)$ we immediately proceed to free motion. Otherwise, the mass stops. If this happens for $\varphi\in [t_1,t_2)$ the mass stops until $\varphi=t_2$ and then switches to free motion. Note, that progression cannot be stopped while $t\in [t_0,t_1)$. If it is stopped on $[t_3,t_0)$ the mass waits the next instant $t_0+2\pi k$ and then starts moving according to Eq.\,(4). In this case, we have $\theta_0=t_0$ and $x_1=0$ in Eq.\,(6). Consequently,
$$\dot x (t)=-a(\cos t-\cos t_0)+ 2c(t-t_0).\eqno (7)$$
This function increases (and, therefore, cannot vanish) until $t=t_1$. However, since $b\ll q$, we may say that $\cos t_3>\cos t_0$ and the right hand side of Eq.\, (7) is negative for $t=t_3$. So, the motion stops somewhere at $[t_1,t_3)$ and then proceeds to the free flight regime. The stop regime may be finished by a transition to a free motion at $t=t_2+2\pi k$ or by transition to the motion with a progression at $t=t_0+2\pi k$. It is impossible neither on $(t_0,t_1)$ nor on $(t_2,t_3)$. See Fig.2 for illustration.

\begin{figure}[!ht]\begin{center}
\includegraphics*[width=2.4in]{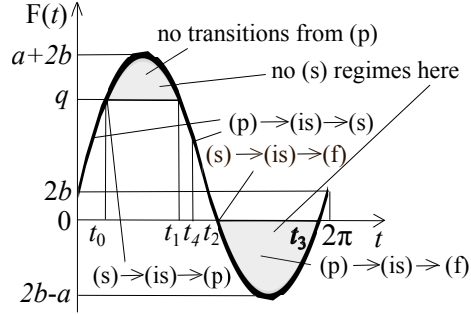}
\end{center}
\caption{Possible regimes of the system and transitions according to the phase.}
\end{figure}

One of the following scenarios must take place for a motion, starting with no-contact regime.

\begin{enumerate}
\item Scenario A: $\mathbf{(f)}\to \mathbf{(fp)}\to \mathbf{(p)}\to \mathbf{(is)} \to \mathbf{(f)}$. This happens if the motion with progression stops at $t\in [t_2+2\pi k,t_3+2\pi k)$.
\item Scenario B: $\mathbf{(f)}\to \mathbf{(fp)}\to \mathbf{(p)}\to \mathbf{(is)} \to \mathbf{(s)}\to \mathbf{(is)}\to \mathbf{(f)}$ -- motion with progression stops at $t\in [t_1+2\pi k,t_2 +2\pi k)$.
\item Scenario C: $$\mathbf{(f)}\to \mathbf{(fp)}\to \mathbf{(p)}\to \mathbf{(is)} \to \mathbf{(s)}\to \mathbf{(is)} \to \mathbf{(p)}\to \mathbf{(is)}\to \mathbf{(s)} \to \mathbf{(is)} \to \mathbf{(f)}.$$
    In this case the first motion with progression stops at one of segments $[t_3+2 \pi k, t_0+2\pi (k+1))$ and a new progression starts at
    $t_0+2\pi (k+1))$ with initial velocity equal to zero. Parameters of the system must be selected so that this second progression stops before
    $t_2+2\pi (k+1)$ otherwise the next scenario is observed
\item Scenario D:
$$\mathbf{(f)} \to \mathbf{(fp)}\to \mathbf{(p)} \to \mathbf{(is)} \to \mathbf{(s)} \to \mathbf{(is)} \to \mathbf{(p)} \to \mathbf{(is)} \to (\ldots).$$
Here $(\ldots)$ implies any sequence of regimes except $\mathbf{(s)}\to \mathbf{(is)}\to \mathbf{(f)}$, corresponding to Scenario C.
\end{enumerate}

Also, there are two degenerate scenarios, corresponding to a zero-velocity stop of the free motion corresponding to $t=t_3+2\pi k$:

\begin{enumerate}
\item Scenario C': $$\mathbf{(f)} \to \mathbf{(is)} \to \mathbf{(s)}\to \mathbf{(is)} \to \mathbf{(p)}\to \mathbf{(is)}\to \mathbf{(s)} \to \mathbf{(is)} \to \mathbf{(f)};$$
\item Scenario D':
$$\mathbf{(f)} \to \mathbf{(is)} \to \mathbf{(s)} \to \mathbf{(is)} \to \mathbf{(p)} \to \mathbf{(is)} \to \mathbf{(\ldots)}.$$
\end{enumerate}

Here we omit all possible instantaneous stops after which the motion returns to the same regime. This does not hurt to equations of motion. However, such stops play an important role since they correspond to discontinuities of stroboscopic mappings. Later on (Section 4) we study them more carefully.

\textbf{Lemma 1}. \emph{There exists a $b_0(q)>0$ such that if $b<b_0(q)$ the following statement is true. Starting progression at the point $t_0$ with an initial velocity equal to zero, the motion must stop at the instant $t_4\in [t_1,t_2]$. Consequently, Scenarios D and D' are impossible for such motions.}

\noindent\textbf{Proof.} If $x(t)$ is a solution of Eq.\,(4) with $\dot x(t_0)=0$, we have
$$\dot x(t_2)=-a\cos t_2+a\cos t_0+2c (t_2-t_0). \eqno (8)$$
In order to prove that this $\dot x(t)$ vanishes somewhere at $[t_0,t_2)$ it suffices to prove that the right hand side of Eq.\, (8) is negative. Instead of this one could prove that
$$a + a \cos t_0 - (q-2b)(\pi-t_0)<0.$$
Here we replaced $t_2$ with $\pi$ in (8) and respected the fact that $2c=2b-q$.
If we demonstrate for a fixed $q$ that
$$a + a \cos t_0 - q(\pi-t_0)<a(1+\cos t_0-\sin t_0(\pi-t_0)) <0.\eqno (9)$$
then there exists a $b_0(q)>0$ such that if $b<b_0(q)$ then the estimate (8) is true.

The second of inequalities (9) is equivalent to the following one: $\cot(t_0/2)<\pi-t_0$ which is true if $t_0>\vartheta_0$ (see Eq.\,(1)) or, equivalently if
$q>\sin \vartheta_0 a$. $\square$

So, wherever and whenever the motion starts, finally, it has a transition to a free motion via an instantaneous stop.

Take an initial instant $\theta\in [t_2,t_3)$ of such transition. There we have $\theta\in [t_2,t_3]$, $x(\theta)=y(\theta)$ (without loss of generality, we may assume that this value is zero) and
$$\dot x(\theta)=\dot y(\theta)=0 \eqno (10)$$
Then the value $\theta$ uniquely defines the farther dynamics.

If $\theta\in [t_2,t_3)$, initially the mass is moving in free regime and, after several transitions, switches to free regime once again. Let $\widehat{T}(\theta)>\theta$ be the first moment of such switching, $T(\theta)=\widehat{T}(\theta) \mod 2\pi$. Both these values are uniquely defined by $\theta$.

So, we may consider the 1D mapping $T: [t_2,t_3)\circlearrowleft$ which is, generally speaking, discontinuous. Considering this mapping only, we lose some information about initial dynamical system, for instance, we do not know any more how the delimiter is shifted.

\medskip

\noindent\textbf{4. Next hit mapping}.

\medskip

Let us introduce an auxiliary mapping $T_1:[t_2,t_3)\to {\mathbb R}$. Let $\theta\in [t_2,t_3)$. For $\theta\in [t_2,t_3)$ we consider a motion with initial conditions (10) and take $\theta_1>\theta$, the first instant when the mass hits the delimiter again. Set $T_1(\theta)=\theta_1$. Note that the image of $T_1$ is an instant, not phase, so it can be greater than $2\pi$.

The value $\theta_1$ corresponds to the first zero of the equation
$$G_1(\theta,t):=b(t-\theta)^2-a\sin t+a\cos \theta(t-\theta)+a\sin\theta=0. \eqno (11)$$
satisfying the condition $t>\theta$. We rewrite this equation in the form
$$\dfrac{b}{a}(t-\theta)^2=\sin t-\cos \theta(t-\theta)-\sin\theta.\eqno (12)$$

The left hand side of Eq.\, (12) is always positive and proportional to the small parameter $b/a$. The right hand side is initially positive and grows faster than the left hand side (both first derivatives vanish for $t=\theta$, but the second derivative of the right hand side is greater since
$\sin \theta\le -b/a<0$).

Geometrically, the right hand side of (12) is the distance between the graph of sine function and the tangent line to it, drawn at $\theta$. If
$\theta\in (3\pi/2,2\pi)$ i.e. $\cos \theta>0$, the graph and its tangent line intersect once again on $(\theta,+\infty)$ and, therefore $\theta_1-\theta<2\pi$. Otherwise, they do not intersect and, for small $b$, there exists a constant
$C>0$ which does not depend on $\theta$ and is such that
$$\theta_1-\theta \ge \dfrac{Ca}{b} \cos\theta.\eqno (13)$$

The mapping $T_1$ is, in general, discontinuous. All possible discontinuities correspond to the case when $\theta_1$ is not a simple zero of (11). In this case derivative $\partial G_1/\partial t$ vanishes for $t=\theta_1$ which means that the following condition is satisfied
$$-a\cos \theta_1+2b(\theta_1-\theta)+a\cos\theta=a\left(2\dfrac{\sin \theta_1 - \sin \theta}{\theta_1-\theta}-\cos\theta - \cos \theta_1\right)=0.\eqno (14)$$
\medskip

\noindent\textbf{5. Points of discontinuity.}

\medskip

\textbf{Lemma 2.} \emph{The intersection of the set of discontinuity points of the mapping $T_1$ with the segment $L_0=[101\pi/100,3\pi/2]$ is finite.}

\noindent\textbf{Proof.} Take $b/a$ so small that $t_2<101\pi/100$. Note that if
$$\theta^{1}<\theta^{2},\qquad \theta^{1,2}\in L_0,\eqno (15)$$
then $T_{1}(\theta^{1})>T_{1}(\theta^{2})$. Indeed, the derivative $\partial G_1/\partial \theta$ of the function $G_1$ defined by Eq.\,(11) equals to
$$-\left(2b+a\sin \theta\right)(t-\theta)$$
and is positive if $\theta\in (t_2,t_3)$ and $t>\theta$. If $\theta^1$ and $\theta^2$ satisfy (15) and $t^1$ is such that $G_1(\theta^1,t^1)=0$ then $G_1(\theta^2,t^1)>0$ and the function $G_1(\theta^2,t)$, negative in a right neighborhood of $t=\theta^2$ must have a zero on $(\theta^2,t^1)$.

So the function $T_1$ is monotonous. Note that if $\theta^1\in L_0$ is a point of discontinuity of $T_1$ and $t^1=T_1(\theta^1)$ then
$$\dfrac{d^2 G_1}{dt^2}(\theta^1,t^1)=2b+a\sin t^1.$$
It follows from Eq.\, (14) that for any discontinuity points $\theta$ of the mapping $T_1$ we must have
$$\cos \theta_1=-\cos \theta +O((\theta_1-\theta)^{-1}).\eqno (16)$$

So, the there exists $\rho>0$ such that the absolute value of second derivative of the function $G_1(\cdot,\theta^1)$ is greater than $\rho$. Consequently, distance between $t^1$ and the next zero of $G(\theta^1,\cdot)$ that is "jump"\ $T_1(\theta^1-0)-T_1(\theta^1+0)$ is greater than a fixed positive value.
This proves that the number of discontinuity points on $L_0$ is finite. $\square$

Grace to Eq.\, (13) and (16) we may take $b/a$ so small that $\theta\in [t_2,5\pi/4]$ then $\theta_1=T(\theta)\in [3\pi/2,t_3]$.

In non-degenerate scenarios ((A)--(C)) at the moment $t=\theta_1$ progression regime starts. The initial velocity of the motion is
$x_1=-a\cos \theta_1+2b(\theta_1-\theta)+a\cos\theta$.
The dynamics of this velocity is described by the formula $\dot x(t)=x_1+2c(t-\theta_1)-a\cos t+a\cos \theta_1$. The progression regime stops as soon as this derivative becomes negative and the next transition $\theta_2$ to free flight or to the stop may be found from equations
$$G_2(\theta,\theta_1,\theta_2):=2b(\theta_1-\theta)+a\cos\theta+2c(\theta_2-\theta_1)-a\cos \theta_2=0. \eqno (17)$$

\textbf{Lemma 3}. \emph{The map $T$ is such that $T(\theta)=\theta_2 \mod 2\pi$ if $\theta_2 \mod 2\pi\in [t_2,t_3)$. Otherwise, a motion with stop has been observed (Scenarios B and C) and $T(\theta)=t_2$.}

\noindent\textbf{Proof.} If $\theta_2\in [t_2,t_3)$ then a motion, starting near the delimiter with the velocity, equal to zero, corresponds to the free regime, so $T=\theta_2$. If $\theta_2\in [t_3,t_0]$, the mass stops until the instant $t_0$ then starts moving in progression regime until $t=t_4$ (see Lemma 1), stops until $t_2$ and proceeds to a free regime. If $\theta_2\in [t_1,t_2)$ the motion stops until $t=t_2$ and also proceeds to the free regime. Since $\theta_2$ cannot belong to $(t_0,t_1)$, we obtain the statement of lemma. $\square$

To finish our proof, we need the following lemma.

\textbf{Lemma 4}. \emph{There exist two disjoint subsegments $J_0$ and $J_1$ of the segment $[t_2,3\pi/2]$ such that
$$T(J_i)\supset [t_2,3\pi/2] \eqno (18)$$
and $T$ is continuous on both segments $J_i$.}

\textbf{Remark.} We may claim without loss of generality that $T(J_i)=[t_2,3\pi/2]$.

\noindent\textbf{Proof.} Let $L_1$ be the arc $[197\pi/100, t_3]$ of the unit circle and $L_2$ be the arc
$$[101\pi/100,51\pi/50].$$
The first arc is correctly defined if the ratio $b/a$ is sufficiently small.

Let $\theta\in L_{1}$. Consider the solution $x(t)$ such that $x(\theta)=y(\theta)=0$, $\dot{x}(\theta)=0$. Then, direct calculations show that for $b=0$, the considered solution starts from free motion and then stops after motion with progression before $t=t_0+2\pi$. So, $\dot x(t_0+2\pi)=0$. Due to continuous dependence of the solution on the parameter $b$, the same is true provided the ratio $b/a$ is sufficiently small.
In this case, as we have already proved $T(\theta)=t_{2}$.

Let $z_{1}<z_{2}<...<z_{n}$ ($n\ge 0$) be discontinuity points of the map $T_{1}$ inside the interval $L_2$. Denote $z_0=101\pi/100$, $z_{n+1}=51\pi/50$.

Suppose that $n<3$. Then there exists $i\in\{0,\ldots, n\}$ such that there is a subsegment ${\cal I}\subset(z_i,z_{i+1})$ of the length not less than $\pi/500$. Let $\theta\in {\cal I}$, $\theta_1=T_1(\theta)$.

It follows from Eq.\,(12) that
$$\theta_1-\theta=\dfrac{a}{b}\left(-\cos \theta+\dfrac{\sin \theta_1-\sin \theta}{\theta_1-\theta}\right).$$
Consequently, if $b/a$ is sufficiently small, the derivative $\partial \theta_1/\partial \theta$ is big on $\cal I$ and values
$$\{\theta_1 \mod 2\pi: \theta\in {\cal I}\}$$
cover $[0, 2\pi)$. Since $T_1$ is continuous on ${\cal I}$, due to (14) we have
$$2\dfrac{\sin \theta_1-\sin\theta}{\theta_1-\theta} -(\cos \theta_1+\cos \theta)\neq 0 \eqno (19)$$
everywhere on $\cal I$. However, due to Eq.\, (13), the maximum of the left hand side of inequality (19) is positive while the minimum is negative. So this inequality cannot hold true everywhere.

So, $n\ge 3$. Then it suffices to prove that $T([z_{i},z_{i+1}))\supset [t_2,3\pi/2]$ for all $i=1,\ldots,n-1$. Note that estimate (16) implies that $T_{1}(z_{i})>3\pi/2$ if $b/a$ is small. Since $\dot{x}(T_{1}(z_{i}))=0$, the corresponding motion proceeds to the free flight immediately after $t=T_{1}(z_{i})$ and, consequently, $T(z_{i})=T_{1}(z_{i})\in L_1$.

On the other hand, $T_{1}(\theta)\longrightarrow T_{1}(T_{1}(z_{i+1}))$ as $\theta\longrightarrow z_{i+1}-0$ (in the limit case, we have a motion, which touches the delimiter with zero velocity). Since $T_{1}(z_{i+1})\in L_{1}$ the corresponding motion is in the stop regime for $t=t_0+2\pi$ and, consequently, $T(\theta)=t_{2}$ for all $\theta$ from a left neighborhood of $z_{i+1}$. This finishes the proof. $\square$

\medskip

\noindent\textbf{6. Infinite set of periodic points.}

\medskip

So, for the mapping $T$ we have obtained two disjoint segments $J_0$ and $J_1$ which are subsets of the arc $[t_2,2\pi/2]$ of the unit circle such that for both $i=0,1$ mappings $T|_{J_i}$ are continuous and $T(J_i)\supset J_0\bigcup J_1$. Let us prove that for any $m\in{\mathbb N}$ the mapping $T$ has a point of the minimal period $m$.

Take a sequence $\{\sigma_k\in \{0,1\}: k\in {\mathbb Z}_+\}$. First of all, we note that there exists a point $p\in J_{\sigma_0}$ such that
$$T^k(p)\in J_{\sigma_k} \eqno (20)$$
for any $k\in {\mathbb N}$.

There exists a segment $J_{\sigma_0\sigma_1}\subset J_{\sigma_0}$ such that $T(J_{\sigma_0\sigma_1})=J_{\sigma_1}$. Then, we may find a segment $J_{\sigma_0\sigma_1\sigma_2}\subset J_{\sigma_0\sigma_1}$ such that $T^2(J_{\sigma_0\sigma_1\sigma_2})=J_{\sigma_2}$. Repeating this procedure, we obtain a nested sequence of segments
$$J_{\sigma_0}\supset J_{\sigma_0\sigma_1} \supset J_{\sigma_0\sigma_1\sigma_2}\supset \ldots$$
The corresponding intersection is non-empty and, consequently, contains a desired point $p$ which may be non-unique.

Fix a number $m$ and consider the sequence $\sigma$, obtained by infinite repetition of a finite sequence $0,\ldots, 0, 1$ of the length $m$. Let
$I=J_{0\ldots 0 1}$ (see above). Then $I \subset T^m(I)$ and the mapping $T^m$ is continuous on this segment. Applying Weierstrass principle to the continuous function $T^m(x)-x$ on the segment I, we obtain a periodic point. Clearly, this point cannot be one of a lower period.

\medskip

\noindent\textbf{7. Symbolical patterns.}

\medskip

Inclusion (18) implies more than just existence of infinite set of periodic points. Here we demonstrate that the reduction of the map $T$ to the union $J_0\bigcup J_1$ is topologically semi-conjugated to the shift of one-side sequences of boolean values.

Let
$$\Sigma=\{\sigma=\{\sigma_k\in\{0,1\}:k\in{\mathbb Z}_+\}\}.$$
Introduce the metrics $d$ on the set $\Sigma$ by the formula
$$d(\sigma,\varsigma)=\sum_{k=0}^\infty 2^{-k} |\sigma_k-\varsigma_k|.$$

Let $J$ be the set of all points $p\in J_0\bigcup J_1$ such that $T^k(p)\in J_0\bigcup J_1$ for all $k\in {\mathbb N}$. Clearly, this set is non-empty and compact. For any $p\in J$ we may introduce the sequence $H(p)=\{\sigma_k\}\in \Sigma$ where values $\sigma_k$ are uniquely defined by Eq.\, (20).

The map $H:J\to \Sigma$ is continuous since all iterations of the map $T|_J$ are continuous. If $S$ is the left shift of sequences of $\Sigma$, defined by the formula
$$S(\sigma)=\varsigma \quad \Leftrightarrow \quad \varsigma_k=\sigma_{k+1}\quad \forall k\in {\mathbb Z}_+$$
the maps $T|_J$ and $S$ are semi-conjugated:
$$H\circ T|_J=S\circ H.$$

\medskip

\noindent\textbf{8. Superstable fixed points.}

\medskip

Here we discuss sufficient conditions for existence of a stable fixed point of the map $T$. Namely, this will be the point $t_2$.

Take $T_1(t_2)>t_2$ i.e. the the first zero of the equation
$$G_1(t_2,t)=0$$
which is a particular case of Eq.\,(11) and $\theta_2(t_2)>\theta_1$ is the first zero of the equation
$$G_2(t_2,T_1(t_2),t)=0$$
which is a particular case of Eq.\,(17).

Let one of inclusions
$$\theta_2(t_2)\mod 2\pi\in (t_1,t_2) \eqno (21)$$
or
$$\theta_2(t_2)\mod 2\pi\in (t_3,2\pi)\bigcup [0,t_0) \eqno (22)$$
be satisfied. Then, due to Lemma 3 there exists $\varepsilon>0$ such that $T([t_2,t_2+\varepsilon))=\{t_2\}$ which means that the point $t_2$ is a \emph{superstable} fixed point of $T$; a neighborhood of this point in $[t_2,t_3)$ is mapped to $t_2$.

Due to Implicit Function Theorem conditions (21) and (22) are robust with respect to small variations of parameters of the considered system if
the following conditions are satisfied:
$$\dfrac{\partial G_1}{\partial \theta}(t_2,\theta)|_{\theta=T_1(t_2)}\neq 0; \qquad
\dfrac{\partial G_1}{\partial \theta}(t_2,T_1(t_2),\theta)|_{\theta=\theta_2(t_2)}\neq 0$$

\textbf{Remark.} If inclusion (21) is true, we do not need to assume that $q>a\sin \vartheta_0$ (see Eq.\, (2)).  Moreover, we do not need the ratio $b/a$ to be small in both cases. We need weaker assumptions
$$a>0, \qquad b>0, \qquad q\in (0, a)\eqno (23)$$
instead.

The obtained superstable periodic solution may coexist with the chaotic invariant set, described in previous sections.

\medskip

\noindent\textbf{9. Conclusion.}

\medskip

Let us formulate principle results of the paper as theorems. Recall that the external force $F(t)$ equals $a\sin t+2b$, $t_2$ and $t_3$ are zeros of $F(t)$ inside $[0,2\pi)$, $q\in (0,a+2b)$ is the maximal value of the dry friction force, $t_0$ and $t_1$ are zeros of $F(t)-q$.

\textbf{Theorem 1}. \emph{For all $a$ and $q$, satisfying inequalities $(2)$ there exists a $b_0=b_0(a,q)>0$ such that for all $b\in (0,b_0)$ the mechanical system, described by equations (3), (4) and (5) is chaotic in the following sense. The phase of transition to free motion uniquely defines the phase of the next transition. This defines a discontinuous mapping $T$ from the segment $[t_2,t_3)$ into itself. There exist two disjoint segments $J_0$ and $J_1$ of the segment $[t_2,3\pi/2]$ such that $T(J_i)\supset [t_2,3\pi/2]$ and $T$ is continuous on both segments $J_i$. Particularly, there exists an infinite set $P$ of periodic points of the mapping $T$. Minimal periods of points of $P$ are unbounded. Moreover, there exists a compact subset $J\subset J_0\bigcup J_1$ such that the map $T|_J$ is continuously semi-conjugated with one-sided symbolic dynamics.}

\textbf{Theorem 2}. \emph{Let inequalities (2) and (22) or inequalities (21) and (23) be satisfied. Then the point $t_2$ is superstable i.e. there exists $\varepsilon>0$ such that $T([t_2,t_2+\varepsilon)=\{t_2\}$.}

\textbf{Remark}. We always assumed that $b<a/2$. If this is not true, any motion eventually does not have free regimes. Then we cannot define map $T$. However, in this case the dynamics of the considered system is very simple. If $q\ge a+2b$ any solution eventually resides in the stop regime. If $q\le 2b$ there exists an instant when an "eternal"\ motion with progression starts. If $q\in (2b-a,a+2b)$ any motion eventually behaves in one of following ways (depending on parameters of the system but not on initial conditions): either it always moves with progression or motion with progression starting at $t_0+2\pi k$ ($k\in {\mathbb Z}$) stops somewhere between $t_1+2\pi k$ and $t_0+2\pi (k+1)$, then motion with progression starts at $t_0+2\pi (k+1)$ and so on.

\medskip

\noindent\textbf{10. Discussion and plans.}

\medskip

First of all, let us note that the obtained chaotic dynamics is robust. Of course, we cannot say anything about stability of points of the set $P$. Every particular point of this set may appear or disappear if we slightly change parameters $a$, $b$ and $q$. The cardinality of the set $P$ may be continuum for some values of parameters, while this set is countable for other values. Neither, the method we offer does not specify the topological structure of the set $P$ and one of its closure.

However, we can select a family of segments $J_0$ and $J_1$ from the statement of Theorem 1 so that boundary points of these segments locally continuously depend on parameters $a$, $b$ and $q$. This can be proved similarly to Lemma 2. The fact of presence of the one dimensional turbulence in the considered system is robust. The same is true for the fact of existence of an infinite set $P$. Moreover, for all fixed values $a$, $b$ and $q$, satisfying inequalities (2), there exists an $\varepsilon>0$ such that for any $C^2$ function
$$G(t,x,\dot x):S^1\times {\mathbb R}^2\to {\mathbb R}$$
such that
$$|G(t,x,\dot x)|<\varepsilon, \qquad \left|\dfrac{\partial G}{\partial (t,x,\dot x)} \right|<\varepsilon$$
for all $t$, $x$ and $\dot x$ an analog of Theorem 1 is true for the system, where equation (3) is replaced with
$$\ddot x=F(t)+G(t,x,\dot x); \quad \dot y=0,$$
equation (4) is replaced with
$$\ddot x=F(t)+G(t,x,\dot x)-q; \quad \dot y=0.$$
and equation (5) is the same.

Particularly, the presence of the considered chaotic dynamics must be observed in simulations and experiments. However, in this paper, we are not going to study the general case. We just offer a method how a non-classical chaos may be found. We plan to use this methods for more general systems with a dry friction
(see \cite{kb18} as an example) and provide for these "real life systems"\ theoretical results accompanied with simulations and experimental data.

\bigskip

\noindent\textbf{Acknowledgements.} This second author was supported in part by Russian Foundation for Basic Researches, grant 12-01-00275-a, by Centre for Research and by FEDER funds through COMPETE – Operational Programme Factors of Competitiveness ("Programa Operacional Factores de Competitividade") and by Portuguese funds through the Center for Research and Development in Mathematics and Applications and the Portuguese Foundation for Science and Technology ("FCT $-$ Funda\c{c}\~{a}o para a Ci\^{e}ncia e a Tecnologia"), within project PEst-C/MAT/UI4106/2011 with COMPETE number FCOMP-01-0124-FEDER-022690. Authors are grateful to Prof. Ron Chen for his precious remarks.

\end{document}